%% file: root.tex
\documentclass[letterpaper, 10 pt, conference]{ieeeconf}  

\IEEEoverridecommandlockouts                             
\usepackage{mathtools,amssymb,commath}
\usepackage{amsmath}
\usepackage{graphicx,import}
\usepackage{hyperref}

\usepackage{mathtools,amsthm}

\hypersetup{colorlinks=true}

\title{\LARGE \bf
Stability 
Analysis and 
Design of \\ Momentum-based 
Controllers  for  Humanoid Robots
}

\author{Gabriele Nava, Francesco Romano, Francesco Nori and Daniele Pucci$^{1}$
\thanks{*This paper was supported by the FP7 EU project KoroiBot (No. 611909 ICT 2013.10 Cognitive Systems and Robotics)}
\thanks{$^{1}$ The authors are with the iCub Facility department, Istituto Italiano di Tecnologia,
        Via Morego 30, Genoa, Italy
        {\tt\small name.surname@iit.it}}%
}

\include{definitions}

\begin{document}

\maketitle
\thispagestyle{empty}
\pagestyle{empty}

\begin{abstract}
    Envisioned applications for humanoid robots call for the design of  balancing and walking controllers. 
    While promising results have been recently achieved, robust and reliable controllers are still a challenge for the control community dealing with humanoid robotics. 
    Momentum-based strategies have proven their effectiveness for controlling humanoids balancing, but the stability analysis of these controllers is still missing. 
    The contribution of this paper is twofold. First, we numerically show that the application of state-of-the-art momentum-based control strategies may lead to 
    unstable zero dynamics. Secondly, we propose simple modifications to the control architecture that avoid instabilities at the zero-dynamics level. Asymptotic 
    stability of the closed loop system is shown by means of a Lyapunov analysis on the linearized system's joint space. 
    The theoretical results are validated with both simulations and experiments on the iCub humanoid robot.
\end{abstract}

\import{tex/}{intro}

\import{tex/}{background}
\import{tex/}{evidence}
\import{tex/}{stability}
\import{tex/}{experiments}
\import{tex/}{conclusions}

\addtolength{\textheight}{-1.2cm}     

\import{tex/}{appendix}

\bibliographystyle{IEEEtran}
\bibliography{IEEEabrv,Bibliography}

\end{document}

%% file: definitions.tex
\newtheorem{Remark}{\bf{Remark}}
\newtheorem{lemma}{\bf{Lemma}}
\newtheorem{assumption}{\bf{Assumption}}

\DeclareMathOperator*{\argmin}{argmin}

\makeatletter
\newsavebox\myboxA
\newsavebox\myboxB
\newlength\mylenA
\newcommand*\xoverline[2][0.75]{%
    \sbox{\myboxA}{$\m@th#2$}%
    \setbox\myboxB\null
    \ht\myboxB=\ht\myboxA%
    \dp\myboxB=\dp\myboxA%
    \wd\myboxB=#1\wd\myboxA
    \sbox\myboxB{$\m@th\overline{\copy\myboxB}$}
    \setlength\mylenA{\the\wd\myboxA}
    \addtolength\mylenA{-\the\wd\myboxB}%
    \ifdim\wd\myboxB<\wd\myboxA%
       \rlap{\hskip 0.5\mylenA\usebox\myboxB}{\usebox\myboxA}%
    \else
        \hskip -0.5\mylenA\rlap{\usebox\myboxA}{\hskip 0.5\mylenA\usebox\myboxB}%
    \fi}
\makeatother

%% file: tex/intro.tex
\section{INTRODUCTION}
\label{sec:intro}

Humanoid robotics is doubtless an emerging field of engineering. One of the reasons accounting for this interest is the need of conceiving systems that can operate 
in places where humans are forbidden to access.
To this purpose, the scientific community has paid much attention to  endowing robots with two main capabilities: locomotion and manipulation. 
At the present day, the results shown at the DARPA Robotics Challenge are promising, but still far from being fully reliable in real applications. 
Stability issues of both low and high level controllers for balancing and  walking were among the main contributors to  robot failures. 
A common high-level control strategy adopted during the competition was that of regulating the robot's momentum, which is usually referred to 
as \emph{momentum-based control}.
This paper presents  numerical evidence that  momentum-based controllers may lead to unstable zero dynamics and proposes a modification to this control scheme
that ensures asymptotic stability.

Balancing controllers for humanoid robots have long attracted the attention of the robotic community~\cite{Caux98,Hirai98}.
Kinematic and dynamic controllers have been common approaches for ensuring a \emph{stable} 
robot behavior for years~\cite{Hyon2007},\cite{Quiang2000}. 
The common denominator of these strategies is  considering the robot attached to ground, which allows one for the application of  classical algorithms 
developed for fixed-based manipulators.

At the modeling level, the emergence of floating-base formalisms for characterizing the dynamics of multi-body systems has loosened 
the assumption of having a robot link attached to ground~\cite{Featherstone2007}. 
At the control level, instead, one of the major complexities when dealing with floating base systems comes from the robot's underactuation. 
In fact, the underactuation forbids the full feedback linearization of the underlying system~\cite{Acosta05}. 
The lack of actuation is usually circumvented by means of rigid contacts between the robot and the environment, but this requires close attention to the 
forces the robot exerts at the contact locations.
If not regulated appropriately, uncontrolled contact forces may break the contact, and the robot control becomes
critical~\cite{Ott2011},\cite{Wensing2013}.

Contact forces/torques, which act on the system as external wrenches, have a direct impact on the rate-of-change of the robot's momentum. 
Indeed, the time derivative of the robot's momentum equals the net wrench acting on the system. 
Furthermore, since the robot's linear momentum can be expressed in terms of the center-of-mass velocity, controlling the robot's momentum is particularly 
tempting for ensuring both robot and contact stability (see, e.g., ~\cite{Frontiers2015} for proper definition of contact stability). 

Several \emph{momentum-based control} strategies have been implemented in real 
applications~\cite{Stephens2010},\cite{Herzog2014},\cite{koolen2015design}.
The essence of these strategies is that of controlling the robot's momentum while guaranteeing stable zero-dynamics. 
The latter objective is often achieved by means of a \emph{postural task}, which usually acts \emph{in the null space} of the control of the robot's 
momentum~\cite{Righetti2011a},\cite{Righetti2011},\cite{Nakanishi11}. These two tasks are achieved by monitoring the contact wrenches, which are ensured to 
belong to the associated feasible domains by resorting to quadratic programming (QP) solvers \cite{Ott2011},\cite{Wensing2013},\cite{Hopkins2015a}. The control of the zero-dynamics can also be used  to control the robot joint 
configuration~\cite{Nakanishi11}.

The control of the linear momentum is exploited to stabilize a desired position for the robot's center of mass. In contrast, the choice of desired values for the robot angular momentum is still unclear \cite{Wieber2006}.
Hence, control strategies that neglect the control of the robot's angular momentum have also been implemented~\cite{Mingxing2015}. 
Controlling both the robot linear and angular momentum, however, is particularly useful for determining the contact torques, and it has become a common 
torque-controlled strategy for dealing with balancing and walking humanoids. Controlling also the angular momentum results in a more
human-like behavior and better response to perturbations~\cite{Hofmann2009},\cite{Herr467}. 

To the best of the  authors' knowledge, the stability analysis of momentum-based control strategies in the contexts of floating base systems is still missing.
The contribution of this paper goes along this direction by considering a humanoid robot standing on one foot, and is then twofold. 
First, we present numerical evidence that classical momentum-based control strategies may lead to unstable zero dynamics, thus meaning that classical 
postural tasks are not sufficient in these cases. 
The main cause of this instability is the lack of  \emph{orientation} correction terms at the angular momentum level: we show that correction terms of 
the form of angular momentum integrals are sufficient for ensuring asymptotic stability of the closed loop system, which is proved by means of a Lyapunov analysis.
The \emph{postural control}, however, is modified with respect to (w.r.t.) state-of-the-art choices. 
The validity of the presented controller is tested both in simulation and on the humanoid robot iCub.

This paper is organized as follows. Section~\ref{sec:background} introduces the notation,  the system modeling, and also recalls a classical momentum-based control 
strategy. Section~\ref{sec:evidence} presents numerical results showing that momentum based control strategies for humanoid robots may lead to unstable zero dynamics. 
Section~\ref{sec:stability} presents a modification of the momentum based control strategy for which stability and convergence can be proven.
Section~\ref{sec:experiments} discusses the numerical and experimental validation of the proposed approach. Conclusions and perspectives conclude the paper.

%% file: tex/background.tex
\section{BACKGROUND}
\label{sec:background}

\subsection{Notation}
\begin{itemize}
        \item $\mathcal{I}$ denotes an inertial frame, with its $z$ axis pointing against the gravity. The constant $g$ denotes the norm of the gravitational acceleration.
        \item $e_i \in \mathbb{R}^m$ is the canonical vector, consisting of all zeros but the $i$-th component that is equal to one.
        \item Given two orientation frames $A$ and $B$, and two vectors $\prescript{A}{}p, \prescript{B}{}p \in \mathbb{R}^{3} $ expressed in these orientation frames, 
              the rotation matrix $\prescript{A}{}R_B$ is such that $\prescript{A}{}p = \prescript{A}{}R_B  \prescript{B}{}p$. 
        \item Let $S(x) \in \mathbb{R}^{3 \times 3}$ be the skew-symmetric matrix such that $S(x)y = x \times y$, where $\times$ is the cross 
              product operator in $\mathbb{R}^3$. 
        \item Given a function $f(x,y) : \mathbb{R}^n \times \mathbb{R}^m \rightarrow \mathbb{R}^p$, the partial derivative of $f(\cdot)$ w.r.t. the variable $x$ is
              denoted as $\partial_x f(x,y) = \frac{\partial f(x,y)}{\partial x} \in \mathbb{R}^{p \times n}$. 
\end{itemize}

\subsection{Modelling}

It is assumed that the robot is composed of $n+1$ rigid bodies, called links, connected by $n$ joints with one degree of freedom each. We also assume
that the multi-body system is \emph{free floating}, i.e. none of the links has an \emph{a priori} constant pose with respect to the inertial frame. The robot configuration space can then be 
characterized by the \emph{position} 
and the \emph{orientation} of a frame attached to a robot's link, called \emph{base frame} $\mathcal{B}$, and the joint configurations. Thus, the 
configuration space  is defined by $\mathbb{Q} = \mathbb{R}^3 \times SO(3) \times \mathbb{R}^n$.
An element of $\mathbb{Q}$ is then a triplet $q = (\prescript{\mathcal{I}}{}p_{\mathcal{B}},\prescript{\mathcal{I}}{}R_{\mathcal{B}},q_j)$,
where $(\prescript{\mathcal{I}}{}p_{\mathcal{B}},\prescript{\mathcal{I}}{}R_{\mathcal{B}})$ denotes the origin  and orientation of the \emph{base frame} expressed 
in the inertial frame, and $q_j$ denotes the \emph{joint angles}. It is possible to define an operation associated with the set $\mathbb{Q}$ such that this set is a 
group. Given two elements $q$ and $\rho$ of the configuration space, the set $\mathbb{Q}$ is a group under the following operation:
$q \cdot \rho = (p_q + p_\rho, R_q R_\rho, q_j + {\rho}_j).$

Furthermore, one easily shows that $\mathbb{Q}$ is  a Lie group. Then, the \emph{velocity} of the multi-body system can be characterized by 
the \emph{algebra} $\mathbb{V}$ of $\mathbb{Q}$ defined by: $\mathbb{V} = \mathbb{R}^3 \times \mathbb{R}^3 \times \mathbb{R}^n$.
An element of $\mathbb{V}$ is then a triplet $\nu = ( ^\mathcal{I}\dot{ p}_{\mathcal{B}},^\mathcal{I}\omega_{\mathcal{B}},\dot{q}_j) = (\text{v}_{\mathcal{B}}, \dot{q}_j)$,
where $^\mathcal{I}\omega_{\mathcal{B}}$ is the angular velocity of the base frame expressed w.r.t. the inertial frame, 
i.e. $^\mathcal{I}\dot{R}_{\mathcal{B}} = S(^\mathcal{I}\omega_{\mathcal{B}})^\mathcal{I}{R}_{\mathcal{B}}$. 

We also assume that the robot is interacting with the environment exchanging $n_c$ distinct wrenches\footnote{As an abuse of notation, we define as \emph{wrench} a quantity that is not the dual of a 
\emph{twist}}. The application of  the Euler-Poincar\'e 
formalism \cite[Ch. 13.5]{Marsden2010} to the multi-body system  yields the following equations of motion: 
\begin{align}
    \label{eq:system}
       {M}(q)\dot{{\nu}} + {C}(q, {\nu}) {\nu} + {G}(q) =  B \tau + \sum_{k = 1}^{n_c} {J}^\top_{\mathcal{C}_k} f_k
\end{align}
where ${M} \in \mathbb{R}^{n+6 \times n+6}$ is the mass matrix, ${C} \in \mathbb{R}^{n+6 \times n+6}$ is the Coriolis matrix, ${G} \in \mathbb{R}^{n+6}$ is the gravity 
term, $B = (0_{n\times 6} , 1_n)^\top$ is a selector matrix, $\tau \in \mathbb{R}^{n}$ is a vector representing the  actuation joint torques, 
and $f_k \in \mathbb{R}^{6}$ denotes the $k$-th external wrench applied by the environment on the robot. 
We assume that the application point 
of the external wrench is associated with a frame $\mathcal{C}_k$, attached to the link on which the wrench acts, and has its $z$ axis pointing in the direction
of the normal of the contact plane. Then,  the external wrench $f_k$ is expressed in a frame whose orientation is that of the inertial 
frame $\mathcal{I}$, and whose origin is that of $\mathcal{C}_k$, i.e. the application point of the external wrench $f_k$. 
The Jacobian $ {J}_{\mathcal{C}_k}=  {J}_{\mathcal{C}_k}(q)$ is the map between the robot's velocity ${\nu}$ and the linear and angular 
velocity $ ^\mathcal{I}\text{v}_{\mathcal{C}_k} := (^\mathcal{I}\dot{ p}_{\mathcal{C}_k},^\mathcal{I}\omega_{\mathcal{C}_k})$ of the frame $\mathcal{C}_k$, i.e.
$^\mathcal{I}\text{v}_{\mathcal{C}_k} = {J}_{\mathcal{C}_k}(q) {\nu}.$ The Jacobian has the following structure: 
\begin{IEEEeqnarray}{RCLRLL}
\label{eqn:jacobianSingle}
{J}_{\mathcal{C}_k}(q) &=& \begin{bmatrix} {J}_{\mathcal{C}_k}^b(q) & {J}_{\mathcal{C}_k}^j(q)\end{bmatrix} &\in& \mathbb{R}^{6\times n+6}, \IEEEyessubnumber \\ 
 {J}_{\mathcal{C}_k}^b(q) &=& 
 \begin{bmatrix}
 1_3 & -S(\prescript{\mathcal{I}}{}p_{\mathcal{C}_k}-\prescript{\mathcal{I}}{}p_{\mathcal{B}})\\ 
 0_{3\times3} & 1_3 \\ 
 \end{bmatrix} &\in& \mathbb{R}^{6\times6} . \IEEEyessubnumber
\end{IEEEeqnarray}
Lastly, it is assumed that  holonomic constraints act on System~\eqref{eq:system}. These constraints are of the form $ c(q) = 0$,
and may represent, for instance, a frame having a constant pose w.r.t. the inertial frame.
In the case where this frame corresponds to the location at which a contact occurs on a link, we represent the holonomic constraint as
${J}_{\mathcal{C}_k}(q) \nu{ =} 0.$

\subsection{Block-Diagonalization of the Mass Matrix}
\label{sec:decoupling}
This section recalls a new expression of the  equations of motion~\eqref{eq:system}. In particular, the next lemma presents a change of coordinates in the
state space $(q,{\nu})$ that transforms the system dynamics~\eqref{eq:system} into a new form where the mass matrix is block diagonal, thus decoupling joint
and base frame accelerations. The obtained equations of motion are then used in the remaining of the paper.

\begin{lemma}
\label{lemma:silvioTranfs}
The proof is given in \cite{traversaro2016}. Consider the equations of motion given by~\eqref{eq:system} and the mass matrix partitioned as following
\[
M = \begin{bmatrix}
    M_b & M_{bj}\\
    M_{bj}^\top & M_j
\end{bmatrix}
\]
with $M_b \in \mathbb{R}^{6 \times 6}$, $M_{bj} \in \mathbb{R}^{6 \times n}$ and $M_j \in \mathbb{R}^{n \times n}$.
Perform the following  change of state variables:
\begin{align}
q := {q}, ~~ \bar{\nu} := T(q){\nu}, \IEEEyessubnumber
\label{eq:centroidalTransform}
\end{align}
with 
\begin{IEEEeqnarray}{rCL}
 \label{eq:generalStructure}
T &:=& \begin{bmatrix}
        \prescript{c}{}X_{\mathcal{B}} & \prescript{c}{}X_{\mathcal{B}} {M}^{-1}_b {M}_{bj} \\
        0_{n \times 6} & 1_n
    \end{bmatrix}, \IEEEyessubnumber \\
    \prescript{c}{}X_{\mathcal{B}} &:=& \begin{bmatrix}
        1_3 & -S(\prescript{\mathcal{I}}{}p_{c}-\prescript{\mathcal{I}}{}p_{\mathcal{B}}) \\ 
        0_{3 \times 3} & 1_3
    \end{bmatrix}\IEEEyessubnumber
\end{IEEEeqnarray}
where the superscript $c$ denotes the frame with the origin located at the center of mass, and with orientation of $\mathcal{I}$.
Then, the equations of motion with state variables $(q,\xoverline{ \nu})$ can be written in the following form
\begin{align}
    \label{eq:decoupled_system_dynamics}
    \xoverline{M}(q)\dot{\xoverline{\nu}} + \xoverline{C}(q, \xoverline{\nu}) \xoverline{\nu} + \xoverline{G} = B \tau +
    \sum_{i = 1}^{n_c} \xoverline{J}^{\top}_{\mathcal{C}_i} f_i,
\end{align}
with
\begin{IEEEeqnarray}{RCL}
    \xoverline{M}(q) &=& T^{-\top} {M} T^{-1} = \begin{bmatrix}
        \xoverline{M}_b(q) & 0_{6\times n} \\ 0_{n\times 6} & \xoverline{M}_j(q_j)
    \end{bmatrix}, \IEEEyessubnumber \label{eq:massMatrixStructure}\\
    \xoverline{C}(q,\xoverline{\nu}) &=& T^{-\top}({M}\dot{T}^{-1} + {C}T^{-1}), \IEEEyessubnumber \label{eq:coriolisMatrixStructure} \\
    \xoverline{G} &=& T^{-\top}{G} = mge_3, \IEEEyessubnumber \label{eq:gravityStructure} \\
    \xoverline{J}_{\mathcal{C}_i}(q) &=& {J}_{\mathcal{C}_i}(q) T^{-1} = \begin{bmatrix} \xoverline{J}_{\mathcal{C}_i}^b(q) & 
    \xoverline{J}_{\mathcal{C}_i}^j(q_j)\end{bmatrix},\IEEEyessubnumber \label{eq:jacobianStructure}
\end{IEEEeqnarray}
 \begin{align*}
    \xoverline{M}_b(q) &= \begin{bmatrix}
        m 1_3 & 0_{3\times3} \\ 0_{3\times3} & I(q)
    \end{bmatrix},
 \xoverline{J}_{\mathcal{C}_i}^b(q) {=}
 \begin{bmatrix}
 1_3 & {-}S(p_{\mathcal{C}_i}{-}\prescript{\mathcal{I}}{}p_{c})  \\ 
 0_{3\times3} & 1_3 \\ 
 \end{bmatrix}
\end{align*}
where $m$ is the mass of the robot and $I$ is the inertia matrix computed with respect to the center of mass, with the orientation of $\mathcal{I}$.
\end{lemma}

The above lemma points out that the mass matrix of the transformed system~\eqref{eq:decoupled_system_dynamics} is block diagonal, i.e. the transformed base acceleration 
is independent from the joint acceleration. More precisely, the transformed robot's velocity $\xoverline{\nu}$ is given by  
$\xoverline{\nu} =
\begin{pmatrix}
^\mathcal{I}\dot{p}_c^\top &
^\mathcal{I}\omega_c^\top &
\dot{q}_j^\top
\end{pmatrix}^\top$
where $^\mathcal{I}\dot{p}_c$ is the velocity of the center-of-mass of the robot, and
$^\mathcal{I}\omega_c$ is the so-called \emph{average angular velocity}\footnote{The term $^\mathcal{I}{\omega}_c$ is also known as 
the \emph{locked angular velocity}~\cite{marsden1993reduced}.}\cite{jellinek1989separation},\cite{essen1993average},\cite{orin2013}.
Hence,
Eq.~\eqref{eq:decoupled_system_dynamics} unifies what the specialized robotic literature usually presents with two sets of equations: the equations of 
motion of the free floating system and the \emph{centroidal dynamics}\footnote{In the specialized literature, the terms \emph{centroidal dynamics} are used 
to indicate the rate of change of the robot's momentum expressed at the center-of-mass, which then equals the summation of all external wrenches acting on 
the multi-body system \cite{orin2013}.}  when expressed in terms of the \emph{average angular velocity}. For the sake of correctness, let us remark that defining 
the \emph{average angular velocity} as the angular velocity of the multi-body system is not theoretically sound. In fact, the existence of a rotation 
matrix $R(q)\in SO(3)$ such that  $\dot{R}(q)R^\top(q) = S(^\mathcal{I}\omega_c)$, i.e. the integrability of $^\mathcal{I}\omega_c$,  is still an open issue.
  
Observe also that the gravity term $\xoverline{G}$ is constant and influences the acceleration of the center-of-mass only. This is a direct consequence of 
\eqref{eq:centroidalTransform}-\eqref{eq:generalStructure} and of the property that $G(q) = Mge_3$, with $e_3 \in \mathbb{R}^{n+6}$.


\begin{Remark}
In the sequel, we assume that the equations of motion are given by 
~\eqref{eq:decoupled_system_dynamics}, 
i.e. the mass matrix is block diagonal. As an abuse of notation but for the sake of simplicity, we hereafter drop the overline notation. 
\end{Remark}
 
\subsection{A classical momentum-based control strategy}
This section recalls a classical momentum-based control strategy when implemented as a two-layer stack-of-task. We assume that the objective 
is the control of the robot momentum and the stability of the zero dynamics.

Recall that the configuration space of the robot evolves in a group of \emph{dimension}\footnote{With group dimension we here mean the dimension of the 
associated algebra $\mathbb{V}$.} $n+6$. Hence, besides pathological cases, when  the system is subject to a set 
of holonomic constraints of dimension~$k$, the configuration space shrinks into a space of dimension $n+6-k$. The stability analysis of the constrained system may 
then require to determine the minimum set of coordinates that characterize the evolution of the constrained system. This operation is, in general, far from obvious
because of the topology of the group $\mathbb{Q}$.

\noindent Now, in the case the holonomic constraint is of the form
$T(q) = \text{constant}$, 
with $T(q) \in \mathbb{R}^3 \times SO(3)$, i.e. a robot link has a constant position-and-orientation w.r.t. the inertial frame, one gets rid of the topology 
related problems of $\mathbb{Q}$ by relating  the base frame $\mathcal{B}$ and the constrained frame. In this case, the minimum set of coordinates belongs 
to $\mathbb{R}^n$ and can be chosen as the joint variables $q_j$.
In light of the above, we make the following assumption.
\begin{assumption}
\label{ass:feet_fixed}
	Only one frame associated with a robot link has a constant position-and-orientation with respect to the inertial frame.
\end{assumption}

Without loss of generality, it is  assumed that the only constrained frame is that between the environment and one of the robot's feet.  Consequently, one has:
\begin{align}
    \label{eq:extForces}
     \sum_{k = 1}^{n_c} {J}^\top_{\mathcal{C}_k} f_k ={J}^\top(q) f,
\end{align}
where $J(q) \in \mathbb{R}^{6\times n+6}$ is the Jacobian of a frame attached to the foot's sole in contact with the environment, 
and $f \in \mathbb{R}^{6}$ the contact wrench.
Differentiating the kinematic constraint 
\begin{IEEEeqnarray}{RCL}
	\label{JNuEqualZero}
    J(q) \nu = \begin{bmatrix}
       J_b &  J_j
    \end{bmatrix}\nu = 0
\end{IEEEeqnarray}
 associated with the contact, yield
\begin{equation}
    \label{eq:constraints_acc}
    \begin{bmatrix}
       J_b &  J_j
    \end{bmatrix} \begin{bmatrix}
                     \dot{\text{v}}_{\mathcal{B}}\\ 
                     \ddot{q}_j
                 \end{bmatrix} + \begin{bmatrix}
                                     \dot{J}_b & \dot{J}_j
                                 \end{bmatrix} \begin{bmatrix}
                                                 \text{v}_{\mathcal{B}}\\
                                                 \dot{q}_j
                                             \end{bmatrix} = 0.
\end{equation}

\subsubsection{Momentum control}
Thanks to the results presented in Lemma~\ref{lemma:silvioTranfs}, the robot's momentum $H \in \mathbb{R}^6$ is given by
    $H = M_b \text{v}_\mathcal{B}$.
The rate-of-change of the robot momentum equals the net external wrench acting on the robot, which in the present case reduces to the contact
wrench $f$ plus the gravity wrench. To control the robot momentum, it is assumed that the contact wrench $f$ can be chosen at will.
Note that given the particular form of~\eqref{eq:decoupled_system_dynamics}, the first six rows correspond to the dynamics of the robot's momentum, i.e.
\begin{IEEEeqnarray}{RCL}
	\label{hDot}
    \frac{\dif }{\dif t}(M_b {\text{v}_\mathcal{B}}) &=& J_b^\top f - mge_3 = \dot{H}(f)
\end{IEEEeqnarray}
where ${H}:=(H_L,H_\omega)$, with 
$H_L, H_\omega \in \mathbb{R}^3$ linear and angular momentum, 
respectively.
The control objective can then be  defined as the stabilization of a desired robot momentum $H^d$. Let us define the momentum error as follows
$\tilde{H} = H - H^d$.
The  control input $f$ in Eq.~\eqref{hDot} is chosen so as
\begin{IEEEeqnarray}{RCL}
	\label{hDotDes}
    \dot{H}(f) &=& 
    \dot{H}^* := \dot{H}^d - K_p \tilde{H} - K_i 
	I_{\tilde{H}}    
    \IEEEeqnarraynumspace  \IEEEyessubnumber \\
    \dot{I}_{\tilde{H}} &=& \tilde{H} \IEEEyessubnumber  \label{eq:IhTilde} 
\end{IEEEeqnarray}
where $K_p, K_i \in \mathbb{R}^{6\times 6}$ are two symmetric, positive definite matrices. It is important to note that a classical choice for the 
matrix $K_i$ consists in~\cite{Herzog2014},\cite{Lee2010}:
\begin{IEEEeqnarray}{RCL}
   K_i =
   \begin{pmatrix}
   	K^L_i & 0_{3 \times 3} \\
   	0_{3 \times 3} & 0_{3 \times 3}
   \end{pmatrix},
	\label{Ki}
\end{IEEEeqnarray}
i.e. the integral correction term at the angular momentum level is equal to zero, while the positive definite matrix $K^L_i~\in~\mathbb{R}^{3\times 3}$ is used for 
tuning the tracking of a desired center-of-mass position when the initial conditions of the integral in~\eqref{hDotDes} are properly set.

Assumption \ref{ass:feet_fixed} implies that the contact wrench satisfying  
Eq.~\eqref{hDotDes} can be chosen as
\begin{equation}
    \label{eq:forces}
    f = J_b^{-\top} \left( 
    \dot{H}^*
    + mg e_3\right) .
\end{equation}

Now, to determine the control torques that instantaneously realize the contact force given 
by~\eqref{eq:forces}, we use the dynamic equations~\eqref{eq:decoupled_system_dynamics} along with the constraints~\eqref{eq:constraints_acc}, which yield
\begin{equation}
    \label{eq:torques}
    \tau = \Lambda^\dagger (JM^{-1}(h - J^\top f) - \dot{J}\nu) + N_\Lambda \tau_0
\end{equation}
where $\Lambda = J_j M_j^{-1} \in \mathbb{R}^{6\times n}$,  $N_\Lambda \in \mathbb{R}^{n\times n}$ is the nullspace projector of $\Lambda$, $h \in \mathbb{R}^{n+6}$
is the vector containing both the Coriolis and gravity terms and $\tau_0 \in \mathbb{R}^n$ is a free variable.

\subsubsection{Stability of the zero dynamics}
The stability of the zero dynamics is usually attempted by means of a so called ``postural task'', 
which exploits the free variable $\tau_0$.
A classical state-of-the-art choice for this postural task consists in:
    $\tau_0= h_j - J_j^\top f -  K^j_{p}(q_j-q_{j}^d) -K^j_{d}\dot{q}_j$,
where  
$h_j - J_j^\top f$ compensates for the nonlinear effect and the external wrenches acting on the joint space of the system.
Hence, the (desired) input torques that are in charge of stabilizing both a desired robot momentum $H_d$ and the associated zero dynamics  are given by
\begin{IEEEeqnarray}{RCL}
	\label{inputTorquesSOT}
	\tau &=& \Lambda^\dagger (JM^{-1}(h - J^\top f) - \dot{J}\nu) + N_\Lambda \tau_0
	\IEEEyessubnumber \IEEEeqnarraynumspace \\
	f &=& J_b^{-\top} \left( 
     \dot{H}^*
	+ mg e_3\right) 
	\IEEEyessubnumber \IEEEeqnarraynumspace \\
    \tau_0 &=& h_j - J_j^\top f -  K^j_{p}(q_j-q_{j}^d) -K^j_{d}\dot{q}_j
    \IEEEyessubnumber \IEEEeqnarraynumspace
\end{IEEEeqnarray}
We present below numerical results showing that~\eqref{inputTorquesSOT} with $K_i$ as~\eqref{Ki} may lead to unstable zero dynamics.

%% file: tex/evidence.tex
\section{NUMERICAL EVIDENCE OF UNSTABLE ZERO DYNAMICS}
\label{sec:evidence}

\subsection{Simulation Environments}

Two different simulation setups have been exploited to perform the numerical validation. In both cases, we simulate the humanoid robot iCub with 23 DoFs \cite{Metta20101125}.

\subsubsection{Custom setup}
\label{sec:num_int_description}
It is in charge of integrating the dynamics~\eqref{eq:decoupled_system_dynamics} when it is subject to the constraint~\eqref{eq:constraints_acc}.
We parametrize $SO(3)$ by means of a quaternion representation 
$\mathcal{Q} \in \mathbb{R}^4 $. 
The resulting  state space system, which is integrated through time, is then:
$
\chi :=(
    p_{\mathcal{B}}, \mathcal{Q}, q_j, \dot{p}_{\mathcal{B}}, \omega_{\mathcal{B}}, \dot{q}_j
)
$, and its derivative is given by
$
\dot{\chi} = (
    \dot{p}_{\mathcal{B}},  \dot{\mathcal{Q}}, \dot{q}_j, \dot{\nu})
$

The constraints~\eqref{eq:constraints_acc}, as well as 
$|\mathcal{Q}|=1$, are then enforced during the integration phase, and additional correction terms have been added \cite{Gros2015}.
The system evolution is then obtained by integrating the constrained dynamical system with the numerical integrator MATLAB \emph{ode15s}.

\subsubsection{Gazebo setup}
The Gazebo simulator~\cite{Koenig04} is the other simulation setup used for our tests. 
Of the different physic engines that can be used with Gazebo, we chose the Open Dynamics Engine (ODE).
Differently from the previous simulation environment, Gazebo allows one for more flexibility. Indeed, we only have to specify the model of the robot, and the 
constraints arise naturally while simulating. Furthermore, Gazebo integrates the dynamics with a fixed step semi-implicit Euler integration scheme. 
Another advantage of using Gazebo w.r.t. the custom integration scheme previously presented consists in the ability to test in simulation the same control software
used on the real robot.

\subsection{Unstable Zero Dynamics}

To show that the momentum-based control strategies may lead to unstable zero dynamics, we control the linear momentum of the robot so as to follow a desired center
of mass trajectory, i.e. a sinusoidal reference along the $y$ coordinate with amplitude $0.05\mathrm{m}$ and frequency $0.3\mathrm{Hz}$. The reference $q^d_j$ is set equal to its joint initial value, i.e. $q^d_j = q_j(0)$.

\subsubsection{Tests on the robot balancing on one foot in the custom simulation setup}
we present simulation results obtained by applying  the control laws~\eqref{inputTorquesSOT}  with $K_i$ as in~\eqref{Ki}. It is assumed that 
the left foot is attached to ground, and no other external wrench applies to the robot.

Figure~\ref{fig:unstable_momentum} shows typical simulation results of the convergence to zero of the robot's momentum error, thus meaning that the  wrench \eqref{eq:forces}
ensures the stabilization of ${\tilde{H}}$ towards zero.  Figure~\ref{fig:unstable_matlab_joints}, instead, depicts the joint position error norm   $|q_j - q^d_j|$.
This figure shows  that the norm of the joint angles increases while the robot's momentum is kept equal to zero, which is a classical behavior of unstable zero dynamics.

\subsubsection{Tests when the robot balances on two feet in the Gazebo simulation setup}

Tests on the stability of the zero dynamics have been carried out also in the case the robot stands on both feet. The main difference between the control algorithm
running in this case and that in~\eqref{inputTorquesSOT} resides in the choice of contact forces and in the holonomic constraints acting on the system. 
More precisely, the contact wrench $f$ is now a twelve dimensional vector, and composed of the contact wrenches $f_L,f_R \in \mathbb{R}^6$ between the floor and
left and right feet, respectively. Hence,  $f = \begin{bmatrix} f_L,f_R \end{bmatrix} \in \mathbb{R}^{12}$.  Also, let $J_L,J_R \in \mathbb{R}^{6 \times n+6}$ denote the Jacobian of two frames
associated with the contact locations of the left and right foot, respectively. Then, 
$J = \begin{bmatrix} J^\top_L,J^\top_R \end{bmatrix}^\top \in \mathbb{R}^{12 \times n+6}$.
 By assuming that the contact wrenches can still be used as virtual control input in the dynamics of the robot's momentum $\dot{H}$ in ~\eqref{hDot}, one is 
 left with a six-dimensional redundancy of the contact wrenches to impose ${H}(f) = \dot{H}^*$. We use this redundancy to minimize the joint torques. In the
 language of \emph{Optimization Theory}, the above control objectives can be formulated as follows. 
\begin{IEEEeqnarray}{RCL}
	\IEEEyesnumber
	\label{optTorque}
	f^* &=& \argmin_{f}  |\tau^*(f)| \IEEEyessubnumber  \\
		   &s.t.& \nonumber \\
		   &&Cf < b \IEEEyessubnumber  \label{frictionCones} \\
		   && \dot{H}(f) = \dot{H}^* \IEEEyessubnumber \\
		   &&\tau^*(f) = \argmin_{\tau}  |\tau(f) - \tau_0(f)| 	\label{optPost} 
  \\
		   	&& \quad s.t.  \nonumber \\
		   	&& \quad \quad \ \dot{J}(q,\nu)\nu + J(q)\dot{\nu} = 0
		    \IEEEyessubnumber 	\label{constraintsRigid} \\
		   	&& \quad \quad \ \dot{\nu} = M^{-1}(B\tau+J^\top(q) f {-} h(q,\nu)) \IEEEyessubnumber \\
		   && \quad \quad \ 	\tau_0 = 
		   h_j - J_j^\top f {-}  K^j_{p}(q_j-q_{j}^d) {-}K^j_{d}\dot{q}_j		    \IEEEyessubnumber
		   \yesnumber
\end{IEEEeqnarray}
Note that the additional constraint~\eqref{frictionCones} ensures that the desired contact wrenches $f$ belong to the associated friction cones. Once the 
optimum $f^*$ has been determined, the input torques $\tau$ are obtained by evaluating the expression~\eqref{optPost}, i.e.
\begin{IEEEeqnarray}{RCL}
	\label{optTorqueFinal}
	\tau = \tau^*(f^*)
		   \yesnumber
\end{IEEEeqnarray}

\begin{figure}[ht]
    \begin{minipage}[c]{8.5cm}
    \centering
    \includegraphics[width=.8\columnwidth]{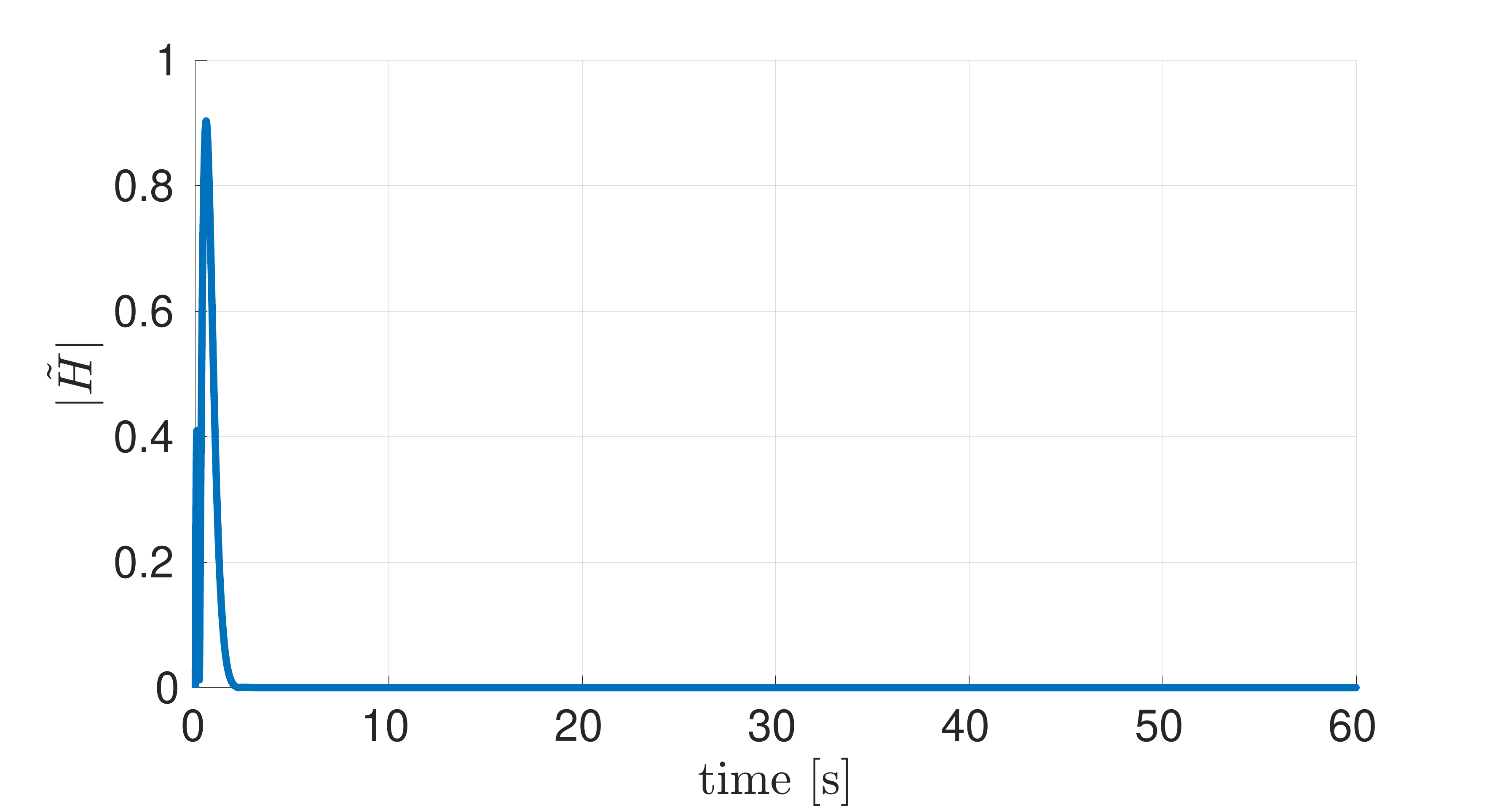}
    \vspace{-0.5cm}
    \caption{Time evolution of the robot's momentum error 
    when standing on one foot and when the control law~\eqref{inputTorquesSOT} is applied. Simulation run with the custom environment.}
    \label{fig:unstable_momentum}
    \end{minipage}
    \begin{minipage}[c]{8.5cm}
    \centering
    \includegraphics[width=.8\columnwidth]{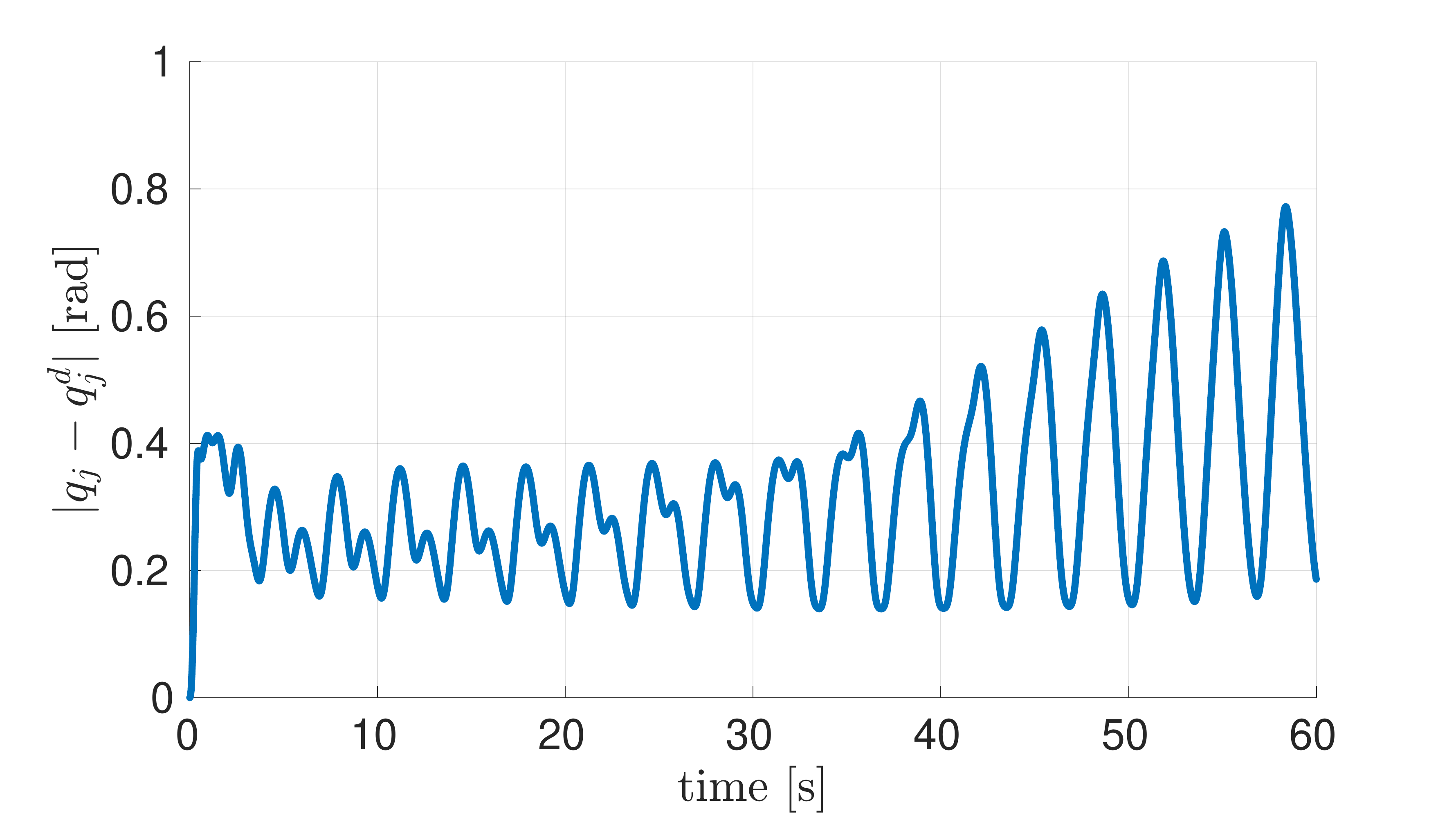}
    \vspace{-0.5cm}
    \caption{Time evolution of the norm of the position error $|q_j - q^d_j|$ 
    when the robot is standing on one foot and when the control law~\eqref{inputTorquesSOT} is applied. Simulation run with the custom environment.}
    \label{fig:unstable_matlab_joints}
    \end{minipage}
    \begin{minipage}[c]{8.5cm}
    \centering
    \includegraphics[width=.8\columnwidth]{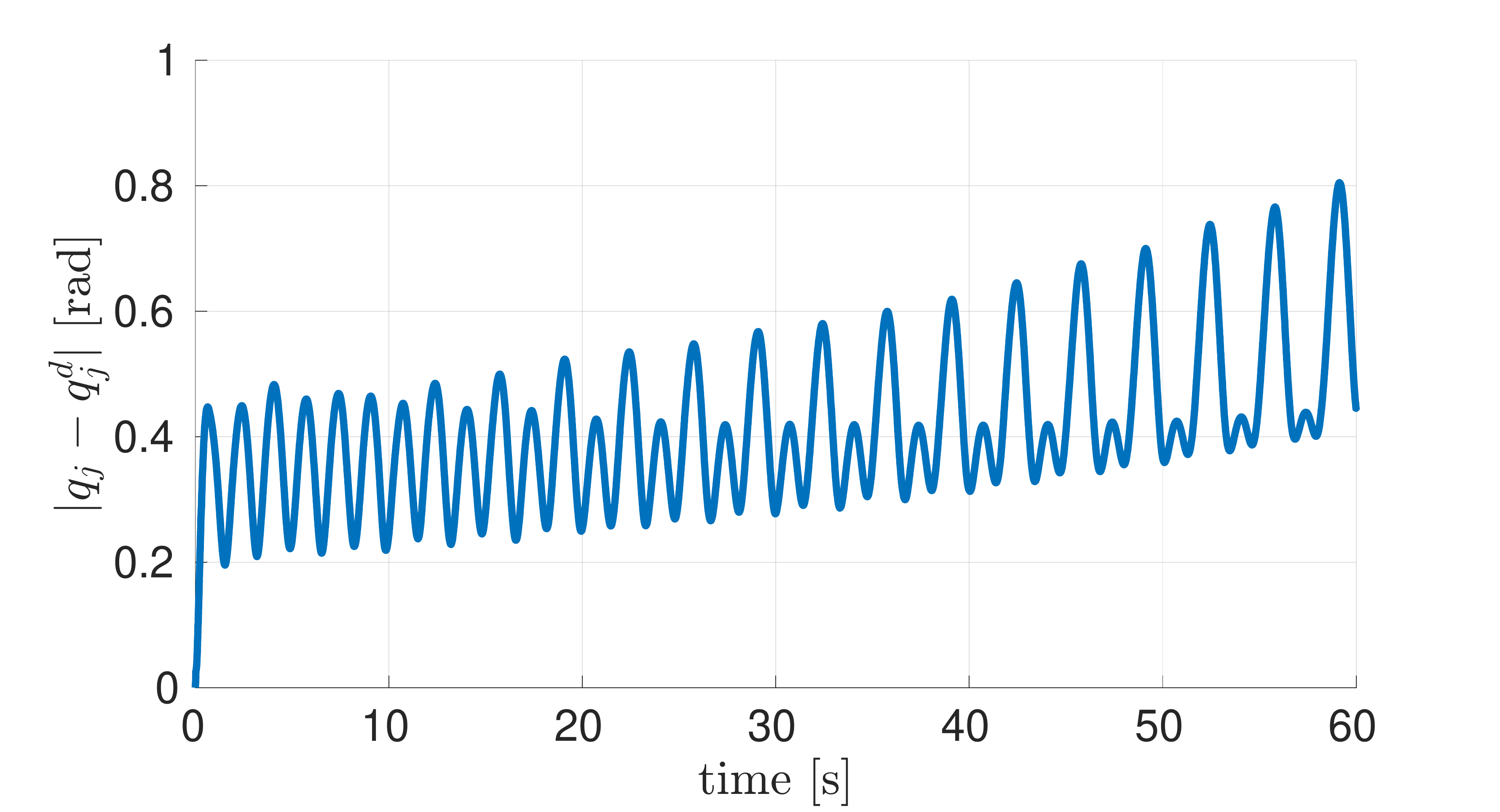}
    \vspace{-0.5cm}
    \caption{Time evolution of the norm of the position error $|q_j - q^d_j|$ 
    when the robot is standing on two feet and when the control law~\eqref{optTorque}-~\eqref{optTorqueFinal} is applied. Simulation run with the Gazebo environment.}
    \label{fig:unstable_simulink_joints}
    \end{minipage}
\end{figure}

Figure~\ref{fig:unstable_simulink_joints} depicts a typical behavior of the joint position error norm  ${|q_j - q^d_j|}$ when the above control algorithm is applied.
It is clear from this figure that the instability of the zero dynamics is observed also in the case where the robot stands on two feet.

%% file: tex/stability.tex
\section{CONTROL DESIGN}
\label{sec:stability}

To circumvent the problems related to the stability of the zero dynamics discussed in the previous section,  
we propose a modification of the control laws~\eqref{inputTorquesSOT} that allows us to show stable zero dynamics of the constrained, closed loop system. The following results exploit the so-called \emph{centroidal-momentum-matrix} $J_G(q) \in \mathbb{R}^{6\times n+6}$, namely the mapping between the robot velocity $\nu$ and the robot's momentum $H$:
    $H = J_G(q)\nu.$
 Now, observe that thanks to the results of Lemma~\ref{lemma:silvioTranfs} one has
${J}_{G}(q) = \begin{bmatrix} M_b & 0_{6 \times n}\end{bmatrix}$.
Then,
observe that when Assumption~\ref{ass:feet_fixed} is satisfied, Eq.~\eqref{JNuEqualZero} allows us to write the robot's momentum  linearly w.r.t. the robot's joint velocity, i.e. 
         $H = \bar{J}_G(q_j)\dot{q}_j,$
where $\bar{J}_G(q_j) \in \mathbb{R}^{6\times n}$ is:
 ˇ\begin{IEEEeqnarray}{RCLRLL}
\label{eqn:jacobian}
\bar{J}_G(q_j)  &:=& 
\begin{bmatrix} {J}_{G}^L(q_j) \\ {J}_{G}^{\omega}(q_j)\end{bmatrix} =
   - M_bJ^{-1}_bJ_j,
\end{IEEEeqnarray} 
and ${J}_{G}^L(q_j),  {J}_{G}^{\omega}(q_j)\in \mathbb{R}^{3\times n}$.
In light of the above, the following result holds.

\begin{lemma}{}
    \label{lemma:stability}
    Assume that Assumption~\ref{ass:feet_fixed} holds, and that the robot possesses more than six degrees of freedom, i.e. $n \geq 6$. In addition, assume 
    also that $H_d = 0$. Let 
	\begin{IEEEeqnarray}{RCL}
		\label{eqPoint}
    		(q_j,\dot{q}_j)&=&(q^d_j,0)
    \end{IEEEeqnarray}
    denote the equilibrium point associated with the constrained, closed loop system and assume that the matrix $\Lambda = J_j(q)M^{-1}_j(q)$ is full row rank in 
    a neighborhood of~\eqref{eqPoint}. 
    Apply the control laws~\eqref{inputTorquesSOT} with  
     \begin{IEEEeqnarray}{RCL}
     	    	\label{modifiedHgDot}
	    \dot{I}_{\tilde{H}}&=& 
	    \begin{bmatrix} {J}_{G}^L(q_j) \\ 
	    {J}_{G}^{\omega}(q^d_j)
	    \end{bmatrix}\dot{q}_j  \IEEEeqnarraynumspace \\
	    	\label{modifiedGains}
		K_i &>& 0,   \quad 
	  K^j_{p} = \xoverline{K}^j_p N_\Lambda M_j,   \quad 
	  K^j_{d}  = \xoverline{K}^j_d N_\Lambda M_j, 
\end{IEEEeqnarray}
where $\xoverline{K}^j_p, \xoverline{K}^j_d \in \mathbb{R}^{n \times n}$ are two constant, positive definite matrices.  
Then, the equilibrium point~\eqref{eqPoint} of the constrained, closed loop dynamics is asymptotically stable.
\end{lemma}
\vspace{0.2cm}
The proof is in Appendix.  Lemma~\ref{lemma:stability} shows that the asymptotic stability of the equilibrium point~\eqref{eqPoint} of the constrained, closed loop dynamics can
be ensured by modifying the integral correction terms, and by modifying the gains of the postural task. As a consequence, the asymptotic stability of the 
equilibrium point~\eqref{eqPoint} implies that the zero dynamics are locally asymptotically stable. 

The fact that the gain matrix $K_i$ must be positive definite conveys the necessity of closing the control loop with \emph{orientation} terms at the angular momentum
level. In fact, some authors intuitively close the angular momentum loop by using the orientation of the  robot's torso~\cite{Ott2011}. 

The proof of Lemma~\ref{lemma:stability} exploits the fact that the minimum coordinates of the robot configuration space when Assumption~\ref{ass:feet_fixed} holds
is given by the joint angles $q_j$. The analysis focuses on the closed loop dynamics of the form
    $\ddot{q}_j = f({q}_j,\dot{q}_j)$
which is then linearized around the equilibrium point~\eqref{eqPoint}. By means of a Lyapunov analysis, one shows that the equilibrium point  is asymptotically stable.
One of the main technical difficulties when linearizing the equation $\ddot{q}_j = f({q}_j,\dot{q}_j)$ comes from the fact that the closed loop dynamics depends on the 
integral of the robot momentum, i.e.
\begin{IEEEeqnarray}{RCL}
	\label{momentumIntegral}
	    {I}_{\tilde{H}}(t) &=& {I}_{\tilde{H}}(0) + \int_0^t
	    \begin{bmatrix} 
	    {J}_{G}^L(q_j(s)) \\ 
	    {J}_{G}^{\omega}(q^d_j)
	    \end{bmatrix}\dot{q}_j(s) \dif s
\end{IEEEeqnarray}
The partial derivative of ${I}_{\tilde{H}}(t)$ w.r.t. the state $(q_j,\dot{q}_j)$  is, in general, not obvious because the matrix ${J}_{G}^L(q_j)$ may not be integrable. Let us observe, however, that the first three rows of  $\dot{I}_{\tilde{H}}(t)$
correspond to the velocity of the center-of-mass times the robot's mass when expressed in terms of the minimal coordinates $q_j$, i.e.
${J}_{G}^L(q_j)\dot{q}_j = m\dot{x}_c$,
with $\dot{x}_c \in \mathbb{R}^3$ the velocity of the robot's center-of-mass. 
Clearly, this means that ${J}_{G}^L(q_j)$ is integrable, and that 
\begin{IEEEeqnarray}{RCL}
	\label{partialDerivatives}
	    \partial_{q_j} {I}_{\tilde{H}} &=& 
	    \begin{bmatrix} 
	    {J}_{G}^L(q_j) \\ 
	    {J}_{G}^{\omega}(q^d_j)
	    \end{bmatrix}   \quad \quad 
	     \partial_{\dot{q}_j} {I}_{\tilde{H}} = 0 \IEEEeqnarraynumspace
\end{IEEEeqnarray}

\begin{Remark}
Lemma~\ref{lemma:stability} suggests that applying the control laws~\eqref{inputTorquesSOT} with the control gains as~\eqref{modifiedGains} can still guarantee stability and convergence of the equilibrium point. First,  observe that the main difference between the variable $I_{\tilde{H}}$ governed by the two expressions~\eqref{eq:IhTilde} and~\eqref{modifiedHgDot} resides only in the last three equations. Then, more importantly, note that the momentum $H$ when Assumption~\ref{ass:feet_fixed} holds can be expressed as follows
	$H = \bar{J}_G(q_j^d)\dot{q}_j + \emph{o} (q_j-q^d_j,\dot{q}_j)$
which implies that 
\begin{IEEEeqnarray}{RCL}
	\label{taylorIntegral}
	 \int_0^t H \dif s &=& \bar{J}_G(q_j^d)(q_j-{q}_j^d) + \int_0^t \emph{o} (q_j-q^d_j,\dot{q}_j) \dif s \nonumber
\end{IEEEeqnarray} 
\end{Remark}
As a consequence, the linear approximations of the integrals $I_{\tilde{H}}$ governed by~\eqref{eq:IhTilde} and~\eqref{modifiedHgDot} coincide when  
\begin{IEEEeqnarray}{RCL}
	\label{assumptionRest}
	 \lim_{ (q_j{-}q^d_j,\dot{q}_j)\to0} \frac{|\int_0^t \emph{o} (q_j-q^d_j,\dot{q}_j)  \dif s|}{|(q_j{-}q^d_j,\dot{q}_j)|} = 0
\end{IEEEeqnarray}
Under the above assumption~\eqref{assumptionRest}, the linear approximation of the control laws~\eqref{inputTorquesSOT} when evaluated with \eqref{eq:IhTilde} and~\eqref{modifiedHgDot} coincide, and stability and convergence of the equilibrium point $(q^d_j,0)$ can still be proven.

%% file: tex/experiments.tex
\section{SIMULATIONS AND EXPERIMENTAL RESULTS}
\label{sec:experiments}

This section shows  simulation and experimental results obtained by applying the control laws~\eqref{inputTorquesSOT}-\eqref{modifiedGains} and~\eqref{optTorque}-\eqref{optTorqueFinal}-\eqref{modifiedGains}. 
To show the improvements of the control modification in Lemma \ref{lemma:stability}, we apply the same reference signal of Section~\ref{sec:evidence}, which revealed unstable zero
dynamics. Hence, the desired linear momentum is chosen so as to follow a sinusoidal reference on the center of mass. Also, control gains are kept equal to those used
for the simulations presented in Section~\ref{sec:evidence}.

\subsection{Simulation results}

\begin{figure}[ht]
   \begin{minipage}[c]{8.5cm}
    \centering
    \includegraphics[width=.8\columnwidth]{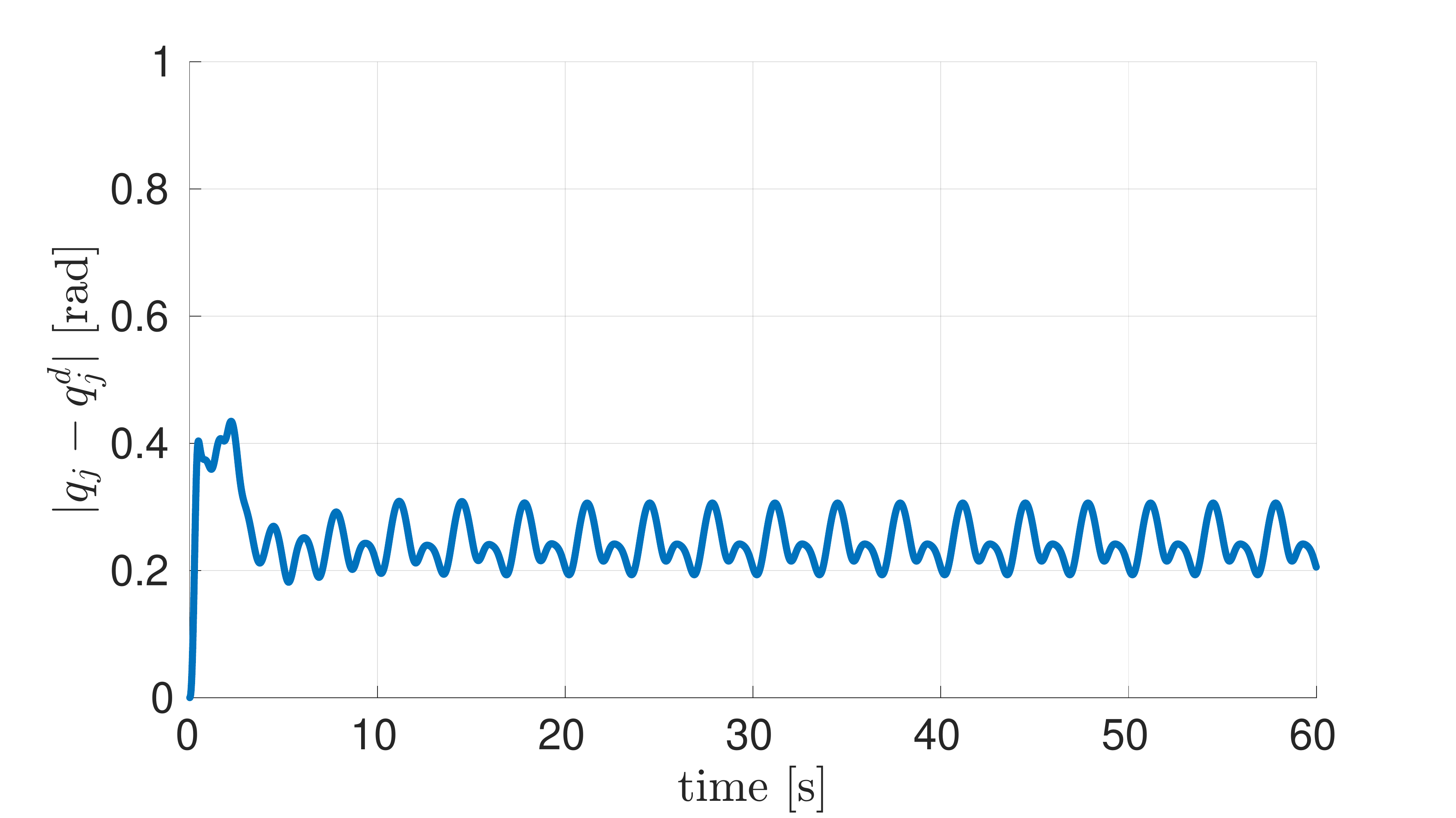}
    \vspace{-0.5cm}
    \caption{Time evolution of the norm of the position error $|q_j - q^d_j|$ 
    when the robot is standing on one foot and when the control law~\eqref{inputTorquesSOT}--\eqref{modifiedGains}--\eqref{modifiedHgDot} is applied. Simulation run with the custom environment.}
    \label{fig:stable_matlab_joints}
    \end{minipage}
   \begin{minipage}[c]{8.5cm}
    \centering
    \includegraphics[width=.8\columnwidth]{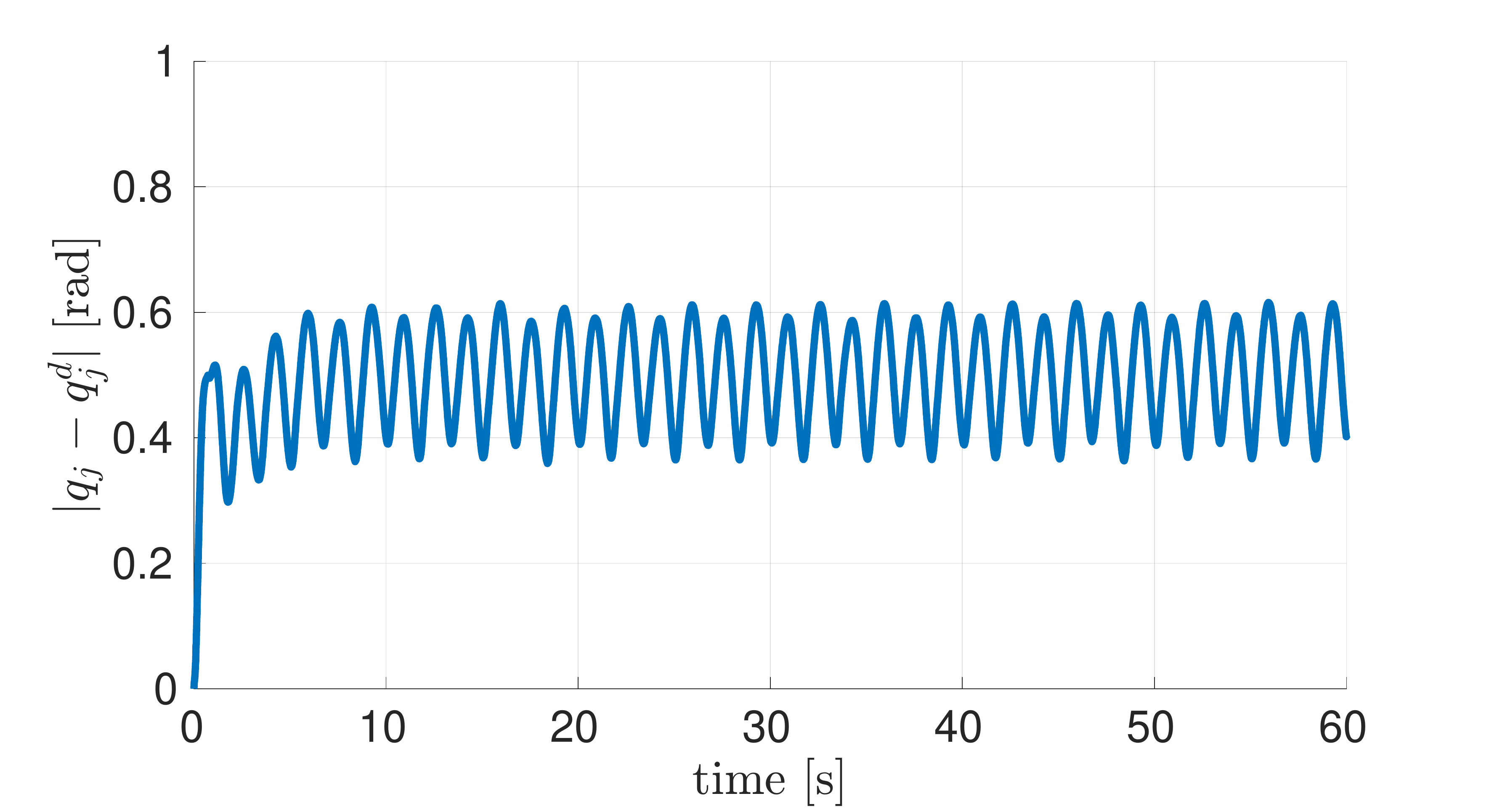}
    \vspace{-0.5cm}
    \caption{Time evolution of the norm of the position error $|q_j - q^d_j|$ 
    when the robot is standing on two feet and when the control law~\eqref{optTorque}--\eqref{optTorqueFinal}--\eqref{modifiedGains}--\eqref{modifiedHgDot} is applied. Simulation run with the Gazebo environment.}
    \label{fig:stable_simulink_joints}
    \end{minipage}
\end{figure}

Figures~\ref{fig:stable_matlab_joints} --~\ref{fig:stable_simulink_joints} show the norm of the joint errors $|q_j-q^d_j|$ when the robot stands on either one or 
two feet, respectively. Experiments on two feet have been performed to verify the robustness of the new control architecture. The simulations are performed both with the custom and Gazebo environment.
As expected, the zero dynamics is now stable, and no divergent behavior of the robot joints is observed.

\subsection{Results on the iCub Humanoid Robot}
We then went one step further and implemented the control algorithm ~\eqref{optTorque}-\eqref{optTorqueFinal} with the modification presented in 
Lemma~\ref{lemma:stability} on the real humanoid robot. The robotic platform used for testing is the iCub humanoid robot~\cite{Metta20101125}.
For the purpose of the proposed control law, iCub is endowed with $23$ degrees of freedom.
A low level torque control loop, running at $1\mathrm{kHz}$, is responsible for stabilizing any desired torque reference signal.

\begin{figure}[ht]
\begin{minipage}[c]{8.5cm}
    \centering
    \includegraphics[width=.8\columnwidth]{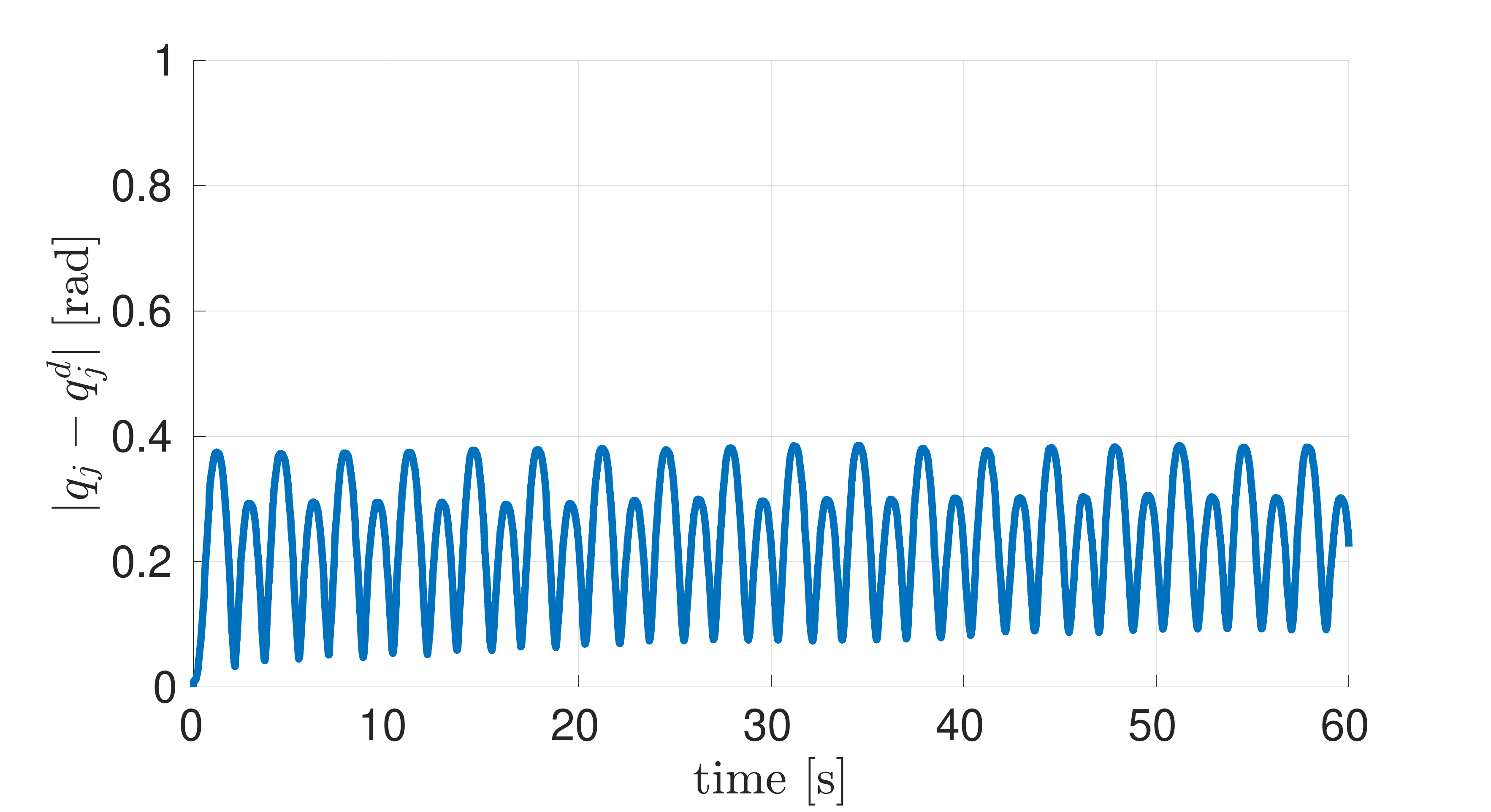}
    \vspace{-0.5cm}
    \caption{Time evolution of the position error norm $|q_j - q^d_j|$ 
    when the robot is standing on two feet and when the laws~\eqref{optTorque}--\eqref{optTorqueFinal}--\eqref{modifiedGains}--\eqref{modifiedHgDot} are applied. Experiment run on the humanoid robot iCub.}
    \label{fig:stable_real_joints}
 \end{minipage}
\begin{minipage}[c]{8.5cm}
    \centering
    \includegraphics[width=.8\columnwidth]{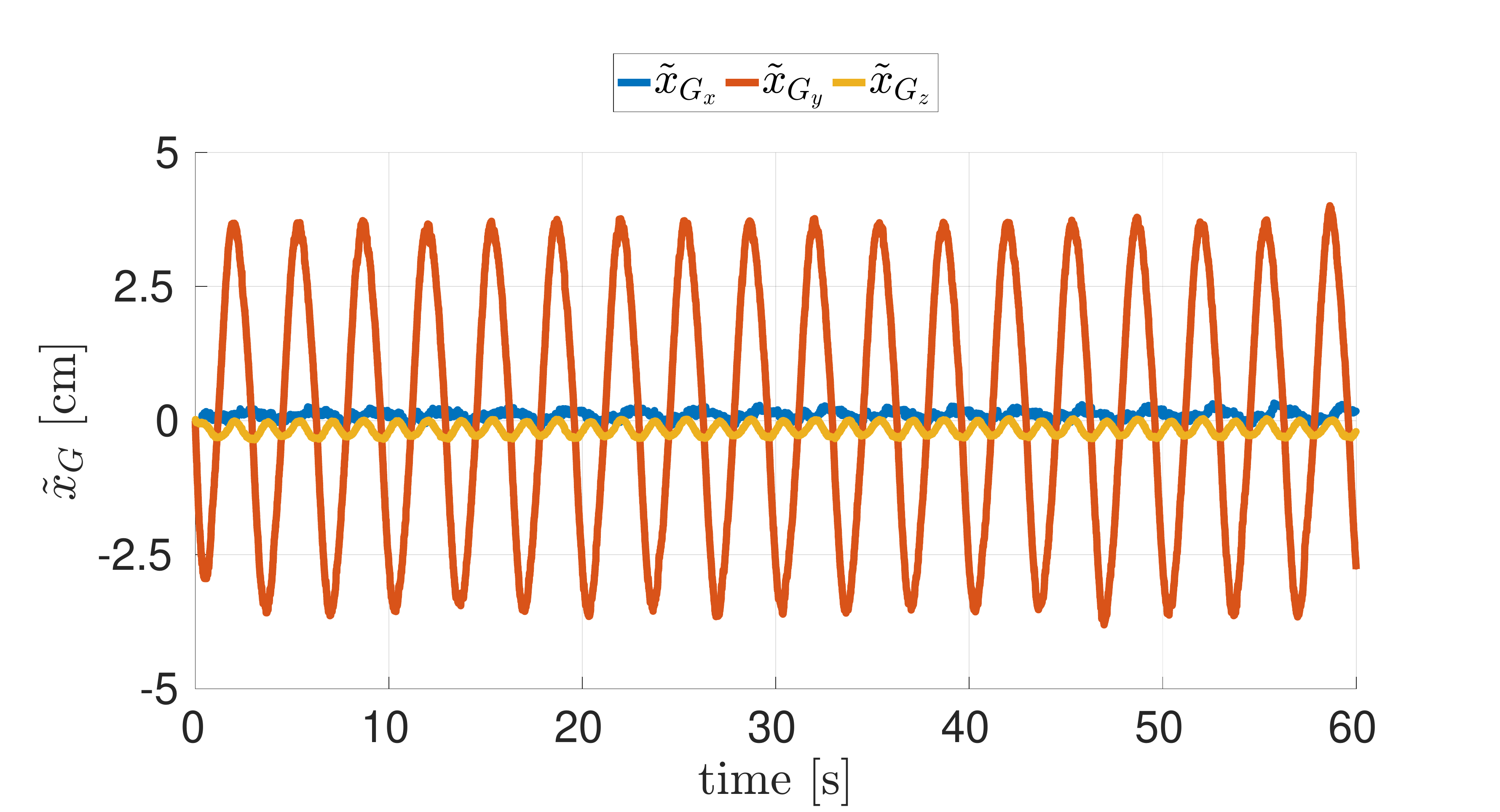}
    \vspace{-0.5cm}
    \caption{Time evolution of robot center-of-mass errors $\tilde{x}_G$
    when the robot stands  on two feet and when the laws~\eqref{optTorque}--\eqref{optTorqueFinal}--\eqref{modifiedGains}--\eqref{modifiedHgDot}  are applied.  Experiment run on the humanoid robot iCub.}
    \label{fig:stable_real_CoMErr}
     \end{minipage}
\end{figure}

Figures~\ref{fig:stable_real_joints} --~\ref{fig:stable_real_CoMErr} show the joint position error $|q_j - q^d_j|$ and the center of mass error. 
Though the center of mass does not converge to the desired value, all signals are bounded, and the control modification presented in Lemma~\ref{lemma:stability} 
does not pose any barrier for the implementation of the control algorithm~\eqref{optTorque}-~\eqref{optTorqueFinal} on a real platform.

%% file: tex/conclusions.tex
\section{CONCLUSIONS AND FUTURE WORKS}
\label{sec:conclusion}
Momentum-based controllers are an efficient control strategy for performing balancing and walking tasks on humanoids.
In this paper, we presented numerical evidence that a stack-of-task approach to this kind of controllers may lead to an instability of the zero dynamics.
In particular, to ensure stability, it is necessary to close the loop with \emph{orientation} terms at the momentum level.
We show a modification of state-of-the-art momentum based control strategies that ensure asymptotic stability, which was shown  by performing Lyapunov analysis 
on the linearization of the closed loop system around the equilibrium point.
Simulation and experimental tests validate the presented analysis.

The stack-of-task approach strongly resembles a cascade of dynamical systems.
It is the authors' opinion that the stability of the whole system can be proved by using the general framework of stability of interconnected systems.
A critical point needed to prove the stability of constrained dynamical system is to define the minimum set  of coordinates identifying its evolution. 
This can be straightforward in case of one contact, but the extension to multiple contacts is not trivial and must be considered carefully.

%% file: tex/appendix.tex
\section*{APPENDIX: proof of Lemma~\ref{lemma:stability}}
\label{proof:lemma_stability}
As described in Section~\ref{sec:stability} the proof is composed of two steps. First, we linearize the constrained closed loop dynamics around the equilibrium 
$(q_j^d, 0)$. Then, by means of Lyapunov analysis, we show that the equilibrium point is asymptotically stable.

\subsubsection{Linearization}
Consider that Assumption~\ref{ass:feet_fixed} holds, and that we apply the control laws~\eqref{inputTorquesSOT}-\eqref{modifiedHgDot} with the gains as~\eqref{modifiedGains}.
The closed loop joint space dynamics of system~\eqref{eq:decoupled_system_dynamics} constrained by~\eqref{eq:constraints_acc} is given by the following equation:
\begin{equation}
    \label{eq:closed_loop_system_joints_vel}
    M_j\ddot{q}_j +h_j -J_j^{\top}f = \tau
\end{equation}
Now, rewrite \eqref{inputTorquesSOT} as follows:
    $\tau = \Lambda^+(\Lambda(h_j-J_j^{\top}f)+J_bM_b^{-1}(h_b-J_b^{\top}f)-\dot{J}\nu)+N_{\Lambda}\tau_0$. 
Therefore, Eq. \eqref{eq:closed_loop_system_joints_vel} can be simplified into:
\begin{equation}
    \label{eq:closed_loop_system_simplified}
    \ddot{q}_j = M_j^{-1}(\Lambda^+(J_bM_b^{-1}(h_b-J_b^{\top}f)-\dot{J}\nu) -N_{\Lambda}u_0)
\end{equation}
where $     u_0   =  K^j_{p}(q_j-q_{j}^d) +K^j_{d}\dot{q}_j$, while $h_b   =  C_b\nu_b + C_{bj}\dot{q}_j + mge_3$, and 
$C_b \in \mathbb{R}^{6\times6}$, $C_{bj} \in\ \mathbb{R}^{6 \times n}$, $C_{jb} \in\ \mathbb{R}^{n \times 6}$,  $C_{j} \in\ \mathbb{R}^{n \times n}$ are obtained from the following partition of 
the Coriolis matrix:
\begin{IEEEeqnarray*} {rCl}
C & = & \begin{bmatrix}
         C_b & C_{bj}\\
         C_{jb} & C_j
        \end{bmatrix}
\end{IEEEeqnarray*}
Substituting \eqref{eq:forces} into \eqref{eq:closed_loop_system_simplified} and grouping together the terms that are linear with respect of joint velocity
yield:
\begin{equation}
    \label{eq:closed_loop_system_final}
    \ddot{q}_j =- M_j^{-1}\left[\Lambda^\dagger (J_b M_b^{-1}\dot{H}^* + \Gamma \dot{q}_j) + N_\Lambda u_0 \right]
\end{equation}
where $    \Gamma = \dot{J}_j - J_bM_b^{-1}C_{bj} + (J_b M_b^{-1} C_b - \dot{J}_b) J_b^{-1} J_j$.\\
Define the state $x$ as
$
    x := \begin{bmatrix}
        x_1^\top & x_2^\top
    \end{bmatrix}^\top = \begin{bmatrix}
        q_j^\top - q^{d\top}_j & \dot{q}_j^\top
    \end{bmatrix}^\top .
$
Since $H^d \equiv 0$, the linearized dynamical system about the equilibrium point $(q_j^d, 0)$ is given by
\begin{equation}
    \label{eq:linear_dynamics}
            \dot{x} = \begin{bmatrix}
                \partial_{q_j} \dot{x}_1 & \partial_{\dot{q}_j} \dot{x}_1 \\
                \partial_{q_j} \dot{x}_2 & \partial_{\dot{q}_j} \dot{x}_2
            \end{bmatrix} x = \begin{bmatrix}
                0_{n \times n} & 1_n \\
                A_1 & A_2
            \end{bmatrix} x
\end{equation}
To find the matrices $A_1,A_2 \in \mathbb{R}^{n\times n}$, one has to evaluate the following partial derivatives
\begin{IEEEeqnarray*}{rCl}
    \partial_{y}\ddot{q}_j &=& - \sum_{i = 1}^6 \partial_y (M_j^{-1} \Lambda^\dagger J_b M_b^{-1} e_i) e_i^\top \dot{H}^* -  M_j^{-1} N_\Lambda  \partial_y u_0  \\
    && -  \sum_{i = 1}^n \partial_y (M_j^{-1} N_\Lambda  e_i) e_i^\top u_0  -  M_j^{-1} \Lambda^\dagger J_b M_b^{-1}  \partial_y\dot{H}^*   
\end{IEEEeqnarray*}
with $y = \{q_j, \dot{q}_j\}$.
Note that $\dot{H}^* = 0$ and $u_0 = 0$ when evaluated at $q_j = q_j^d$ and $\dot{q}_j = 0$.
We thus have to compute only the partial derivatives of $\dot{H}^*$ and $u_0$.
The latter is trivially given by
    $\partial_{q_j} u_0 = \xoverline{K}^j_p N_\Lambda M_j$ and $\partial_{\dot{q}_j} u_0  = \xoverline{K}^j_d N_\Lambda M_j$. The former can be calculated via Eq.~\eqref{partialDerivatives}. In light of the above we obtain the expressions of the matrices in \eqref{eq:linear_dynamics}:
%
\begin{IEEEeqnarray*}{rCl}
    A_1 &=& - M_j^{-1} \Lambda^\dagger J_b M_b^{-1} K_i M_b J_b^{-1}J_j    - M_j^{-1} N_\Lambda  \xoverline{K}^j_p N_\Lambda M_j \\
    A_2 &=& - M_j^{-1} \Lambda^\dagger J_b M_b^{-1} K_p M_b J_b^{-1} J_j  - M_j^{-1} N_\Lambda  \xoverline{K}^j_d N_\Lambda M_j.
\end{IEEEeqnarray*}

\subsubsection{Proof of Asymptotic Stability}
    Consider now the following Lyapunov candidate:
   \begin{IEEEeqnarray}{RCL}
        \label{eq:candidate}
        V(x) &{=}& \frac{1}{2}\left[{x}_1^\top M_j^\top Q_1  M_j {x}_1 {+} x_2^\top M_j^\top Q_2  M_jx_2\right] \nonumber
    \end{IEEEeqnarray}
    where $M_j = M_j(q_j^d)$, and
    \begin{align*}
        Q_1 &:= \Lambda^\top J_b^{-\top}M_b^\top K_i M_b J_b^{-1} \Lambda + N_\Lambda \xoverline{K}_p^j N_\Lambda \\%
        Q_2 &:= \Lambda^\top J_b^{-\top}M_b^\top M_b J_b^{-1} \Lambda + N_\Lambda %
    \end{align*}
    calculated at ${x}_1 = 0$, $x_2 = 0$. $V$ is a properly defined candidate, in fact $V = 0 \iff x = 0$ and is positive definite otherwise.
    Indeed, $Q_1$ can be rewritten in the following way:
    \begin{align*}
        Q_1 = & \begin{bmatrix}
            \Lambda^\top & N_\Lambda
        \end{bmatrix} \begin{bmatrix}
            J_b^{-\top}M_b^\top K_i M_b J_b^{-1}  & 0 \\
            0 & \xoverline{K}_p^j 
        \end{bmatrix} \begin{bmatrix}
            \Lambda\\
            N_\Lambda
        \end{bmatrix}.
    \end{align*}
    and, because $\Lambda$ and $N_\Lambda$ are orthogonal $Q_1$ is positive definite. The same reasoning can be applied to $Q_2$.
    
    We can now consider the time derivative of $V$:
    \begin{IEEEeqnarray*}{rCl}
        \dot{V} &=&  {x}_1^\top M_j^\top Q_1  M_j x_2 + x_2^\top M_j^\top Q_2  M_j\dot{x}_2 \\
            &=&-x_2^\top M_j^\top ( \Lambda^\top J_b^{-\top} M_b^\top K_p M_b J_b^{-1} \Lambda \\
            && + N_\Lambda K^j_d N_\Lambda) M_jx_2 \leq 0.
    \end{IEEEeqnarray*}
    The stability of the equilibrium point $x = 0$ associated with the linear system~\eqref{eq:linear_dynamics} thus follows.
 To prove the asymptotic stability of the equilibrium point $x = 0$, which implies its  asymptotic stability when associated with the nonlinear system~\eqref{eq:closed_loop_system_final},  we have to resort to LaSalle's invariance principle.
    Let us define the invariant set $S:=\{x : \dot{V}(x) = 0 \}$ that implies $S = \{({x}_1, 0)\}$.
    It is easy to verify that the only trajectory starting in $S$ and remains in $S$ is given by ${x}_1 = 0$ thus proving LaSalle's principle.
    As a consequence, the equilibrium point $x = 0 \Rightarrow (q_j, \dot{q}_j) = (q_j^d,0)$ is asymptotically stable.